\documentclass[authoryear,preprint,3p,times,review,12pt]{elsarticle}

\usepackage{xcolor}
\usepackage{subcaption}

\usepackage{amssymb}
\usepackage{amsmath}

\usepackage{amsthm}
\newtheorem{proposition}{Proposition}[section]
\newtheorem{theorem}{Theorem}[section]      
\newtheorem{corollary}[theorem]{Corollary}  

\journal{arXiv}
\begin{document}
\begin{frontmatter}



\title{Strong Consistency of the SIMEX Estimator in Linear Regression with a Conditionally Poisson Covariate} 



{\author[a]{Aijun Yang}
\author[a]{Mary Lesperance}
\author[a]{Farouk S. Nathoo  \corref{cor1}}}
\ead{nathoo@uvic.ca}
\ead{https://nathoogroup.ca}
\cortext[cor1]{Corresponding author: Farouk S. Nathoo}

\affiliation [a]{organization={Department of Mathematics \& Statistics,University of Victoria},
            addressline={3800 Finnerty Rd}, 
            city={Victoria},
            postcode={V8P5C2}, 
            state={BC},
            country={Canada}}

\begin{abstract}
This paper considers estimation for linear regression analysis with covariate measurement error arising from Poisson surrogates. We consider cases where covariates follow a conditional Poisson distribution, capturing non-Gaussian and heteroscedastic error structures. To address this, we extend the simulation extrapolation (SIMEX) algorithm to the conditional Poisson setting (POI-SIMEX), enabling robust adjustment in the absence of internal validation data. Theoretical analysis establishes strong consistency of the POI-SIMEX estimator under a linear regression framework.
\end{abstract}




\end{frontmatter}



\section{Introduction} \label{sec1}
The presence of covariate measurement error in regression analysis can distort the estimated coefficients, leading to biased estimators with increased variance \citep{Carroll2006,Fuller1986}. This issue is generally well-studied, but there is limited work considering situations where the observed covariates are discrete and follow a Poisson distribution conditional on the true covariate. This problem is motivated by applications in biostatistics, such as the analysis of tissue
microarrays (TMAs) where marker counts within biological sub-samples serve as error prone proxies for true underlying densities. However, our focus here is on providing a theoretical justification for applying SIMEX in this setting.


The literature on adjustments for covariate measurement error is substantial (see, e.g., \citet{Fuller1986,Carroll2006}, and \citet{GraceYi2021}). In the existing literature, approaches for handling conditional Poisson surrogates in contexts lacking validation data where regression calibration cannot be applied have been largely overlooked. Regression analysis and measurement error modeling approaches often rely on classical errors-in-variables models. These methods attempt to correct bias by explicitly modeling measurement errors in the independent variables, usually assuming Gaussian distributed errors. Regression calibration \citep{Spiegelman1997} is commonly used for bias correction in settings where internal validation data are available.

\citet{Nakamura1990} discusses the corrected score function approach for generalized linear models and shows that the solution to the unbiased score equations will yield a consistent estimator. \citet{Lietal2004} utilize this approach and derive the corrected score function for linear regression with Poisson-distributed covariates. They derive explicit forms for the solutions to the corrected score equation, under the assumption that the true covariate, the mean of the observed covariate, follows a Gamma distribution. While it is suitable for linear regression with Poisson distributed surrogates, the approach of \citet{Lietal2004} does not easily generalize beyond this narrow setting to more complex setups involving censored data or nonlinear regression.

Among existing methods, simulation-extrapolation (SIMEX) requires only weak assumptions about the distribution of the true predictor and has been extended in several directions \citep{Cook1994,Carroll1995}. \citet{Kuchenhoff2006} develop misclassification simulation-extrapolation (MC-SIMEX) to deal with measurement error for categorical variables. \citet{Parveen2021TheNM} extend SIMEX to allow for non-zero mean measurement error.

In certain situations where the variance of the measurement error is unknown, and it is heteroscedastic, \citet{Devanarayan2002} propose an empirical SIMEX estimator where the pseudo errors are generated from random linear contrasts of the observed replicate measurements. A comprehensive review of SIMEX and related developments is presented in \citet{Sevilimedu2022}. 

While the corrected score estimator for Poisson surrogates proposed by \citet{Lietal2004} can be applied in this context, its application is limited to linear regression. In contrast, SIMEX can be applied more broadly, but it has yet to be developed for conditionally Poisson-distributed surrogates. In addition, even in the context of linear regression, simulations have found that SIMEX generally outperforms the corrected score methodology. 

Despite these advances, there is no theoretical justification for SIMEX under a conditionally Poisson measurement error model. The standard SIMEX algorithm assumes homoscedastic, additive error and employs Gaussian pseudo-errors, 
assumptions that are violated in our context of heteroscedastic, discrete surrogates. It is thus unclear whether the method remains valid when applied to this important class of problems.

The primary contribution of this work is to provide the first theoretical justification for a SIMEX-type estimator under a conditionally Poisson measurement error model, establishing conditions for its strong consistency. We formulate the POI-SIMEX algorithm to address this setting and make the following specific contributions: 1) formulation of a conditional Poisson error model for discrete-valued surrogates; 2) extension of the SIMEX algorithm to the Poisson case (POI-SIMEX), accommodating heteroscedastic error variance with only a single replicate; 3) establishing the strong consistency of the POI-SIMEX estimator for linear regression. 

Subsequent sections of the paper proceed as follows. Section 2 focuses on the conditional Poisson error model framework and SIMEX estimation, including model formulation and the proposed POI-SIMEX estimator. In Section 3, we establish the strong consistency of the POI-SIMEX estimator in simple linear regression for a conditionally Poisson-distributed observed covariate. Section 4 presents a numerical example. The paper concludes with a discussion of our results and future work.

\section{Poisson Error Model Framework} \label{sec2}
\subsection{Error model assumptions}
We model $Y$, the response variable, as a function of $X$, the true covariate of interest, where $X$ typically represents a marker density $W_T/A_T$, $W_T$ a count over a larger area $T$ of size $A_T$, cannot be measured directly and instead, the variable $W$, representing the marker count of a selected subregion $A$ is measured. The estimated or observed marker density is $W/A$ and is used as the surrogate for $W_T/A_T$. Other covariates $Z$ are assumed to be measured without error. Marker locations within the selected region are assumed to be distributed according to a homogeneous Poisson process with intensity $X = \frac{W_T}{A_T}$. 
Under the Poisson process model for the marker distribution, we have for subject $i$,
\begin{equation}
\label{eq0}
W_i\lvert X_i \overset{\displaystyle \text{ind}}{\sim} \text{Poisson} (X_i A_i) \ \ i=1,...,N, 
\end{equation}
where $X_i$ is the unobserved true covariate for the $i^{th}$ subject and $W_i$ is the observed surrogate count with area $A_i$. Thus, $E[W_i\lvert X_i]= Var [W_i\lvert X_i]=A_i X_i$, and $E[W_i]=A_i E[X_i]$. We note that the model is Poisson conditional on $X_{i}$, and therefore overdispersion in the observed $W_{i}$ does not preclude the use of this model. 

In what follows, we assume $(Y_i,Z_i,W_i,A_i,X_i)$, corresponding to the response, the covariate without measurement error, a surrogate covariate observed with Poisson noise, the observed subregion area, and the true covariate, respectively, are independent across subjects $i=1,...N$. For simplicity, we illustrate the conditional Poisson measurement error model in the context of a linear model and assume both $X_i$ and $Z_i$ are scalars. We then have
\begin{equation}
Y_i =\beta_0 + \beta_XX_i + \beta_ZZ_i + \epsilon_i, \ \ \epsilon_i \overset{\displaystyle \text{i.i.d}}{\sim} \text{N} (0,\sigma_\epsilon^2) , \ \  W_i\lvert X_i \overset{\displaystyle \text{ind}}{\sim} \text{Poisson} (X_i A_i), 
\label{eq1}
\end{equation}
where $X_i >0$ is the true marker density, and  $(Y_i,Z_i,W_i,A_i)$ $i=1,..., N$ are the observed data.

\subsection{The SIMEX algorithm for Poisson distributed surrogates}
The SIMEX algorithm, introduced by \citet{Cook1994}, is commonly employed to correct the bias from measurement error in regression analysis, operating under certain assumptions regarding the measurement error structure. In practice, the standard SIMEX algorithm assumes the measurement error follows an additive structure with errors distributed as $N(0,\sigma^2)$, where $\sigma^2$ is either known or estimated. The standard SIMEX method involves three main steps: simulation, estimation, and extrapolation. In the simulation step, additional measurement error is artificially introduced into the data. In the estimation step, the model's parameters are estimated repeatedly with varying levels of introduced error. In the final step, the relationship between the introduced error and the estimated parameters is analyzed, and this information is used to extrapolate the estimators back to the scenario with no measurement error. 

In departing with standard SIMEX methodology, the non-Gaussian measurement error structure in our setup is $U_{i} = \frac{W_{i}}{A_{i}} - X_{i}$, where $W_{i}$ is conditionally Poisson distributed, $A_{i}$ is fixed and known, and the distribution of $X_{i}$ is continuous but unspecified. In addition, the measurement error variance is not constant and varies across subjects. Thus our model lies outside the standard assumptions within which SIMEX is applied \citep{Carroll1997}. 

Remarkably, the use of Gaussian pseudo errors in the proposed SIMEX algorithm produces an estimator that exhibits strong consistency (Theorem \ref{thm1}). 

In the simulation step, we simulate 
$b=1,...,B$ new sets of the observed covariate $W_{i,b}(\lambda)=W_{i}+\lambda^{1/2}\hat{\sigma}_i U_{i,b}$ by adding Gaussian error $\lambda^{1/2}\hat{\sigma_i}U_{i,b}$, where $\hat{\sigma}_i^2$ is an estimate of $VAR[W_i/A_i-X_i]=\frac{E(X_i)}{A_i}$ and $U_{i,b} \overset{\displaystyle \text{i.i.d}} \sim N(0,1)$. The variance of the total measurement error in $W_{i,b}(\lambda)$ is $\sigma_i^2(1+\lambda)$. The regression coefficient estimates of $W_{i,b}(\lambda)$, $\hat{\beta}_{\lambda, b}$ are computed for each of the $B$ simulated datasets, and the average of these estimates $\frac{1}{B}\sum_{b=1}^{B}\hat{\beta}_{\lambda, b}$ is computed as an approximation to $E_U[\hat{\beta}_{\lambda, b}|{Y,W}]$. These computations are repeated over a grid of $\lambda$ values. After the simulation step, one analyzes how the computed expectations vary as $\lambda$ changes by fitting a functional form, such as linear, quadratic, or nonlinear. The SIMEX estimator is then obtained by extrapolating the resulting curve to the value corresponding to the absence of measurement error, i.e., when $\lambda=-1$, so that $\sigma_i^2(1+\lambda)=0$; thus, the variance of the total measurement error is 0. 

The implementation of SIMEX requires estimation of $\sigma_i^2$, the variance of the measurement error for subject $i$, i=1,\dots,N. A key aspect of our proposed POI-SIMEX and an important issue for a study is consistent estimation of the heteroscedastic measurement error variance from the observed data with only a single replicate.

\begin{proposition} \label{prop1}
If $W_i \lvert X_i \overset{ind}{\sim} Poisson (X_iA_i)$, $X_1,..., X_N$ are i.i.d. with $E[X_i] = m < \infty $, $\bar{V}=\frac{1}{N}\sum_i^{N}V_i$ with $V_i=W_i/A_i$ so that $E[V_i] = E[X_i] = m$ and assume $E[(V_i-m)^{4}] \le a < \infty$,
then $\hat{\sigma}_i^2=\frac{\bar{V}}{A_i}$ is a strongly consistent estimator of $\sigma_i^2=Var[W_i/A_i-X_i],\, i=1,\dots,N$.
\end{proposition}
\begin{proof}
The $V_i=W_i/A_i$ are independent with $E[V_i]=E[E[V_i\lvert X_i]]=E[X_i] < \infty$, and does not depend on $i$ since the $x_i'$s are i.i.d. So, the SLLN applies to $\bar{V}$ with $\bar{V} \overset{\text{a.s.}}{\rightarrow}E[W_i/A_i]$. It suffices to show that $\sigma_i^2=E[W_i/A_i]/A_i$. We have
$\sigma_i^2=Var[W_i/A_i-X_i]=E[Var[W_i/A_i-X_i\lvert X_i]] + Var[E[W_i/A_i-X_i\lvert X_i]]=E[X_i]/A_i=E[W_i/A_i]/A_i$, so that $\hat{\sigma}_i^2 \overset{\text{a.s.}}{\rightarrow} \sigma_i^2, \, i=1,\dots,N$.
\end{proof}

\section{Consistency of the SIMEX estimator with conditionally \\Poisson-distributed surrogates in simple linear regression} \label{sec3}

\citet{Carroll1996} and \citet{Carroll1997} develop the asymptotic theory for SIMEX under the Gaussian and additive measurement error model framework. In this section, we demonstrate the consistency of the SIMEX estimator in simple linear regression models with conditionally Poisson-distributed surrogates. 
Let $Y_i=\beta_0+\beta_1X_i+\epsilon_i$, $\epsilon_i \overset{\text{i.i.d}}{\sim} N(0,\tau^2)$ and $W_i \lvert X_i \overset{\text{ind}}{\sim}$ Poisson $(X_iA_i), i=1,..,N; W_i \bot Y_i \lvert X_i$ with observed data $\{(Y_i,W_i)_{i=1}^{N}\}$. For convenience, we assume $A_i=1$ and consider the SIMEX estimator of $(\beta_0$ $\beta_1)$ defined as:
\begin{align*} 
\hat{\beta}_{SIMEX,\sigma} =\lim_{\lambda\to\ -1}E[\hat{\beta}_{l,U}(\lambda) \lvert \{(Y_i,W_i)_{i=1}^{N}\}],\ \ l=0,1
\end{align*}
\vspace{-1em}
\noindent where 
$$
\hat{\beta}_{1,U}(\lambda)=\frac{\sum\limits_{i=1}^N(Y_i-\bar{Y})(W_{i,U}(\lambda)-\bar{W}_U(\lambda))}{\sum\limits_{i=1}^N(W_{i,U}(\lambda)-\bar{W}_U(\lambda))^2} \text{, }
\hat{\beta}_{0,U}(\lambda) =\bar{Y} - \hat{\beta}_{1,U}(\lambda) \bar{W}_U(\lambda)
\vspace{1em}
\text {, and } W_{i,U}(\lambda) = W_i+ \lambda^{\frac{1}{2}}\sigma U_i,$$
$$
U_i \overset{iid}{\sim} N(0,1) \text{ with }
 \bar{W}_U(\lambda)=\frac{1}{N} \sum\limits_{i=1}^NW_{i,U}(\lambda),
 \text{ and } \newline \sigma^2=\text{var}[W_i-X_i] \text{ is assumed known.}
$$
\begin{theorem}\label{thm1}
Under the following regularity conditions:

\noindent(R1) $\lim_{N \to\infty}E[\hat{\beta}_{l,U}(\lambda) \lvert \{(Y_i,W_i)_{i=1}^{N}\}]$ converges uniformly in $\lambda$ on a set having $\lambda=-1$ as a limit point, l=0,1.\\
(R2) $\lim_{\lambda \to -1 }E[\hat{\beta}_{l,U}(\lambda) \lvert \{(Y_i,W_i)_{i=1}^{N}\}]$ exists for each $N,l=0,1$ \\
(R3) $\lim_{N \to\infty}E[\hat{\beta}_{l,U}(\lambda) \lvert \{(Y_i,W_i)_{i=1}^{N}\}]=E[\lim_{N \to \infty} \hat{\beta}_{l,U}(\lambda) \lvert \{(Y_i,W_i)_{i=1}^{N}\}]$ \\
(R4) $E[X_iW_i] < \infty, E[X_i]<\infty, E[W_i^2] <\infty $

We have that
$$\lim_{N \to\infty} \lim_{\lambda \to-1} E[\hat{\beta}_{l,U}(\lambda) \lvert \{(Y_i,W_i)_{i=1}^{N}\}]=\beta_l  \text{ a.s.}$$
The SIMEX estimator for simple linear regression with conditional Poisson measurement error is therefore strongly consistent when $\sigma^2$ is known.
\end{theorem}

\begin{proof}
Given (R1 and R2), we have from Theorem 7.11 \citep{Rudin1976} that 
\begin{equation}
\begin{split}
\lim_{N \to\infty} \lim_{\lambda \to -1}E[\hat{\beta}_{l,U}(\lambda) \lvert \{(Y_i,W_i)_{i=1}^{N}\}]&=\lim_{\lambda \to -1}\lim_{N \to \infty} E[\hat{\beta}_{l,U}(\lambda) \lvert \{(Y_i,W_i)_{i=1}^{N}\}]\\
&=\lim_{\lambda \to -1} E[\lim_{N \to \infty} \hat{\beta}_{l,U}(\lambda) \lvert \{(Y_i,W_i)_{i=1}^{N}\}] ,
\end{split}
\label{eq3}
\end{equation}
the later equality being from (R3).

Consider now the limit, $\lim_{N \to \infty} \hat{\beta}_{l,U}(\lambda)$ for l=1.
\begin{align*} 
 \hat{\beta}_{1,U}(\lambda)&=\frac{\sum_{i=1}^N(Y_i-\bar{Y})(W_{i,U}(\lambda)-\bar{W}_U(\lambda))}{\sum_{i=1}^N(W_{i,U}(\lambda)-\bar{W}_U(\lambda))^2}
=\frac{\sum_{i=1}^N(Y_i-\bar{Y})(W_i-\bar{W}+\lambda^{\frac{1}{2}}\sigma U_i-\lambda^{\frac{1}{2}}\sigma \bar{U})}{\sum_{i=1}^N(W_i-\bar{W}+\lambda^{\frac{1}{2}}\sigma U_i-\lambda^{\frac{1}{2}}\sigma \bar{U})^2}\\
\vspace{2em}
&=\frac{\sum_{i=1}^N(\beta_1(X_i-\bar{X})+(\epsilon_i-\bar{\epsilon})(W_i-\bar{W}+\lambda^{\frac{1}{2}}\sigma U_i-\lambda^{\frac{1}{2}}\sigma \bar{U})}{\sum_{i=1}^N(W_i-\bar{W}+\lambda^{\frac{1}{2}}\sigma U_i-\lambda^{\frac{1}{2}}\sigma \bar{U})^2},\\
\vspace{2em}
&=\frac{\left[\beta_1S_{XW}+\beta_1\lambda^{\frac{1}{2}}\sigma S_{XU}+S_{\epsilon W}+\lambda^{\frac{1}{2}}\sigma S_{\epsilon U}\right]}{\left [\sum_{i=1}^N(W_i-\bar{W})^2+2\lambda^{\frac{1}{2}}\sigma \sum_{i=1}^{N}(W_i-\bar{W})(U_i-\bar{U})+\lambda \sigma^2 \sum_{i=1}^{N}(U_i-\bar{U})^2 \right ]},
\end{align*}
where $S_{XW}=\sum_{i=1}^N(X_i-\bar{X})(W_i-\bar{W})$, $S_{XU}=\sum_{i=1}^N(X_i-\bar{X})(U_i-\bar{U})$, $S_{\epsilon W}=\sum_{i=1}^N(\epsilon_i-\bar{\epsilon})(W_i-\bar{W})$, \\$S_{\epsilon U}=\sum_{i=1}^N(\epsilon_i-\bar{\epsilon})(U_i-\bar{U})$.\\
Next, multiply the numerator and denominator by $\frac{1}{N}$, and by (R4) we apply the strong law of large numbers (Theorem 5.44 \citep{Ross2006}) as well as the Slutsky theorem (Theorem 6 \citep{Ferguson1996}). Hence,
\begin{align*} 
\hat{\beta}_{1,U}(\lambda) \overset{\text{a.s.}}{\rightarrow} \frac{\left [\beta_1 \text{Cov}[X,W]+ \beta_1 \lambda^{\frac{1}{2}}\sigma \text{Cov}[X,U]+\text{Cov}[\epsilon,W]+\lambda^{\frac{1}{2}}\sigma \text{Cov}[\epsilon,U]\right ]} {\left [\text{Var}[W] + 2\lambda^{\frac{1}{2}}\sigma \text{Cov}[W,U] + \lambda \sigma^2 \text{Var}[U] \right ]}
\end{align*}
We have $\text{Cov}[X,U]=0, \text{Cov}[\epsilon,W]=0, \text{Cov}[\epsilon,U]=0, \text{Cov}[W,U]=0$, \\
so that
$\hat{\beta}_{1,U}(\lambda) \overset{a.s}{\rightarrow} \frac{\beta_1 \text{Cov}[X,W]} {\text{Var}[W]+\lambda \sigma^2}$.\\
Now, under the conditional Poisson measurement error model,\\
given $\text{Var}[W]=E[Var[W\lvert X]]+Var[E[W\lvert X]]
=E[X]+Var[X]$ and 
$Cov[X,W]=E[XW]-E[X]E[W]=E[E[XW|X]]-E[X]E[W]=E[X E[W\lvert X]] -E[X] E[E[W\lvert X]]
=E[X^2]-E^2[X]=Var[X] $, hence,
$\hat{\beta}_{1,U}(\lambda) \overset{a.s}{\rightarrow} \frac{\beta_1 Var[X]}{Var[X] + E[X]+ \lambda \sigma^2}$.\\
Then, note that
$\sigma^2=Var[W_i-X_i]=E[Var[W_i-X_i\lvert X_i]]+Var[E[W_i-X_i\lvert X_i]]=E[X_i]$,
so
$\hat{\beta}_{1,U}(\lambda) \overset{a.s}{\rightarrow}  \frac{\beta_1 Var[X]}{Var[X] + (1+\lambda)E[X]}$.\\
Therefore, we have from Equation \ref{eq3}
\begin{align*} 
\lim_{N\to\infty}\lim_{\lambda \to -1} E[\hat{\beta}_{1,U}(\lambda) \lvert \{Y_i,W_i\}_{i=1}^N]
&=\lim_{\lambda \to -1} E\left [\frac{\beta_1 Var[X]}{Var[X]+(1+\lambda)E[X]} \lvert \{Y_i,W_i\}_{i=1}^\infty \right]\\
&=\lim_{\lambda \to -1} \frac{\beta_1 Var[X]}{Var[X]+(1+\lambda)E[X]}\\
&=\beta_1 \ \ (a.s) 
\end{align*}
Next, consider
\begin{align*} 
\hat{\beta}_{0,U}(\lambda)&=\bar{Y}-\hat{\beta}_{1,U}(\lambda)\bar{W}_U(\lambda)\\
&=\frac{1}{N}\sum_{i=1}^{N}(\beta_0+\beta_1X_i+\epsilon_i)-\hat{\beta}_{1,U}(\lambda)\frac{1}{N}\sum_{i=1}^{N}(W_i+\lambda^{\frac{1}{2}}\sigma U_i)\\
&=\beta_0+\beta_1\bar{X}+\bar{\epsilon} - \hat{\beta}_{1,U}(\lambda) (\bar{W}+\lambda^{\frac{1}{2}}\sigma \bar{U})\\
&=\beta_0+\beta_1 \bar{X}+\bar{\epsilon} - \hat{\beta}_{1,U}(\lambda) \bar{W} -\lambda^{\frac{1}{2}}\sigma \hat{\beta}_{1,U}(\lambda)\bar{U}
\end{align*}
and taking the limit,
\begin{align*} 
\lim_{N\to\infty} \hat{\beta}_{0,U}(\lambda)&=\beta_0+\beta_1 E[X]+E[\epsilon] -\frac{\beta_1 Var[X]E[W]}{Var[X]+(1+\lambda)E[X]} -\frac{\lambda^{\frac{1}{2}} \sigma E[U]\beta_1 Var[X]}{Var[X]+(1+\lambda)E[X]}\\
&=\beta_0+\beta_1 E[X]-\frac{\beta_1 Var[X]E[W]}{Var[X]+(1+\lambda)E[X]} \ \ (a.s.) \ \ (SLLN)
\end{align*}
then by Equation \ref{eq3}, 
\begin{align*} 
\lim_{N\to\infty}\lim_{\lambda \to -1} E[\hat{\beta}_{0,U}(\lambda) \lvert \{Y_i,W_i\}_{i=1}^{\infty}=\beta_0+\beta_1 E[X]-\beta_1 E[W]
=\beta_0 \ \ (a.s), 
\end{align*}
since $E[W]=E[E[W|X]]=E[X]$.\\
so that $\hat{\beta}_{SIMEX,\sigma}$ is strongly consistent for $\beta_0, \beta_1$ when $\sigma$ is known.
\end{proof}

\begin{corollary} \label{cor:varEst}
Under the conditions of Theorem (\ref{thm1}) and Proposition (\ref{prop1}) with $A_i=1$, we have that $\hat{\beta}_{SIMEX,\hat{\sigma}}$ with $\hat{\sigma}^2=\bar{W}$, is strongly consistent. The proof proceeds in the same way as that of Theorem (\ref{thm1}), with $\bar{W}$ replacing $\sigma^2$ and further application of Slutsky's theorem when taking the limit.
\end{corollary}

\section{Numerical Example}
\label{sec4}
\begin{figure}[ht]  
  \centering
  \begin{subfigure}[b]{0.47\textwidth}
    \includegraphics[width=\textwidth]{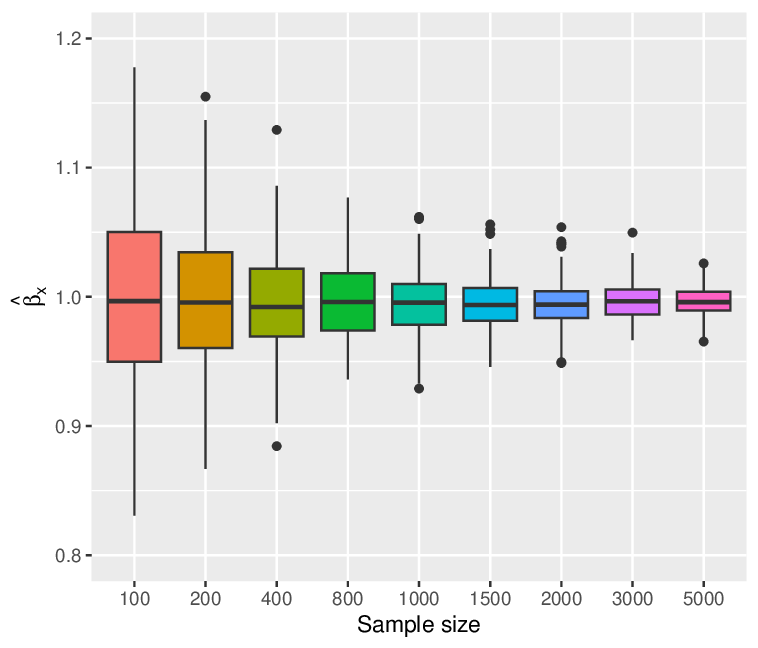}
    \caption{Use true measurement error variance}
    \label{fig:sub1}
  \end{subfigure}
  \hfill
  \begin{subfigure}[b]{0.47\textwidth}
    \includegraphics[width=\textwidth]{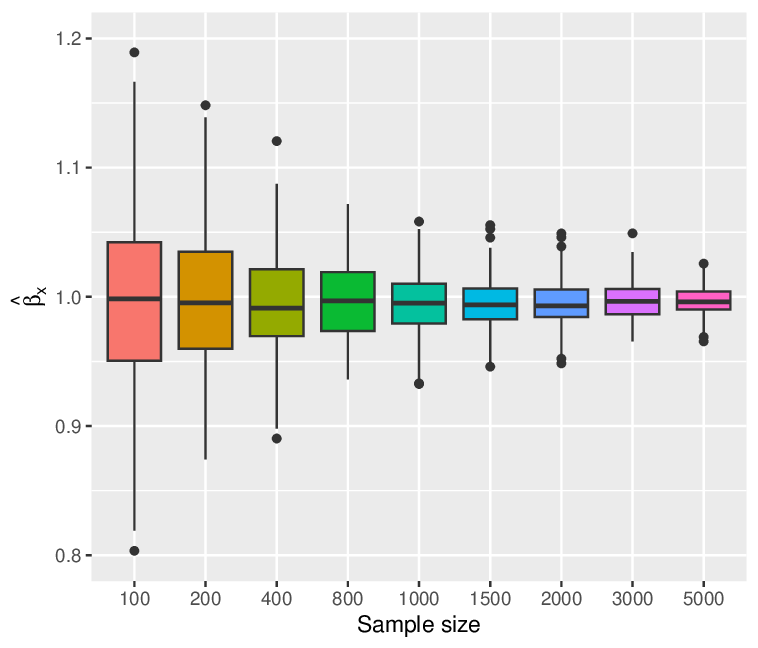}
    \caption{Use estimated measurement error variance}
    \label{fig:sub2}
  \end{subfigure}
  \caption{Boxplots of POI-SIMEX estimates for the regression coefficient of the error-prone covariate across 200 simulations, stratified by sample size, using the true (a) and estimated (b) measurement error variance.}
  \label{fig:fig1a_b}
\end{figure}

We present a numerical example to demonstrate the large sample behavior of the POI-SIMEX estimator. The data are simulated assuming the linear regression model (\ref{eq1}) with the observed count $W$ in an area $A$ distributed as Poisson with mean $XA$ where $X$ is the true covariate. The surrogate for the true covariate is $W/A$ where $W \sim Poisson(XA)$. For simplicity, the area $A$ is set to 1 in the simulations. Covariate $Z$ is measured without error and is generated from a uniform distribution between 0.5 and 9. The coefficients associated with the intercept, covariate with measurement error and without measurement error, respectively $\beta_0$, $\beta_X$, $\beta_Z$, are set to 2, 1, and 0.5, following \citet{Lietal2004}. The sample size N is set to be $100, 200,400,800, 1000,1500,2000,3000,5000$. The standard deviation of the error ($\sigma_\epsilon$) in the regression model is set to 5. The true covariate $X_i$’s are generated independently from a Gamma distribution with shape parameter, $a=1$, and scale parameter, $b=10$. For each sample size, 200 simulated datasets are generated. Each dataset yields one POI-SIMEX estimate of $\beta_x$, and the boxplots in Figure~\ref{fig:fig1a_b} summarize these estimates. Figure~\ref{fig:sub1} uses the true variance of the measurement error; whereas, Figure~\ref{fig:sub2} employs the estimated variance obtained following procedure outlined in Proposition \ref{prop1}. The plots demonstrate the narrowing of the samplig distribution of the POI-SIMEX estimator as the sample size increases to $3000$ towards the true value.


\section{Conclusion} \label{sec5}
We developed POI-SIMEX, an adjustment of SIMEX designed to incorporate conditionally Poisson-distributed surrogates, and applied it within the context of linear regression models. For scenarios involving unknown measurement error variance, our method requires only a single replicate to yield a strongly consistent estimate of the measurement error variance, making POI-SIMEX both practical and efficient. A numerical example is used to demonstrate this behavior. Its practical utility and superior finite-sample performance of this estimator compared to existing methods are demonstrated in a separate application-focused paper.

While this paper has initiated the theoretical justification for POI-SIMEX, more general theory including the application of POI-SIMEX for semiparametric survival modeling with the Cox proportional hazard model is part of ongoing work.





\section*{Acknowledgment}
Nathoo and Lesperance acknowledge funding from NSERC (RGPIN$-04044-2020$, RGPIN$-07079-2020$) through the Discovery Grants Program.  Yang acknowledges funding from the University of Victoria and the Pacific Leaders Scholarship Program.

\section*{Data availability}
Only simulated data was used for the research described in the article.

\end{document}